\documentclass[12pt,english]{amsart}
\usepackage[T1]{fontenc}
\usepackage[latin9]{inputenc}
\usepackage{geometry}
\usepackage{color}
\usepackage{babel}
\usepackage{verbatim}
\usepackage{amsthm}
\usepackage{amssymb}
\usepackage[unicode=true]
 {hyperref}
\usepackage{breakurl}

\makeatletter
\numberwithin{equation}{section}
\numberwithin{figure}{section}
\theoremstyle{plain}
\newtheorem{thm}{\protect\theoremname}[section]
  \theoremstyle{plain}
  \newtheorem{fact}[thm]{\protect\factname}
  \theoremstyle{remark}
  \newtheorem{rem}[thm]{\protect\remarkname}
  \theoremstyle{plain}
  \newtheorem{cor}[thm]{\protect\corollaryname}

\@ifundefined{definecolor}
 {\usepackage{color}}{}

\@ifundefined{definecolor}{\usepackage{color}}{}
\usepackage{amsfonts}\usepackage{babel}

\providecommand{\U}[1]{\protect\rule{.1in}{.1in}}

\usepackage{calrsfs}
\DeclareMathAlphabet{\pazocal}{OMS}{zplm}{m}{n}

\usepackage{babel}

\makeatother

  \providecommand{\corollaryname}{Corollary}
  \providecommand{\factname}{Fact}
  \providecommand{\remarkname}{Remark}
\providecommand{\theoremname}{Theorem}

\begin{document}

\title[Complex interpolation of compactness]{{\normalsize Lecture notes on complex interpolation of compactness
- preliminary version}}

\author{{Michael Cwikel}}

\address{Cwikel: Department of Mathematics, Technion - Israel Institute of
Technology, Haifa 32000, Israel }

\email{mcwikel@math.technion.ac.il }

\author{{Richard Rochberg}}

\address{Rochberg: Department of Mathematics, Washington University, St. Louis,
MO, 63130, USA }

\email{rr@math.wustl.edu}

\thanks{The first named author's work was supported by the Technion V.P.R.\ Fund
and by the Fund for Promotion of Research at the Technion. The second
named author's work was supported by the National Science Foundation
under Grant No. 1001488.}
\begin{abstract}
We recall that the fundamental theorem of complex interpolation is
the

\textbf{Boundedness Theorem:} If, for $j=0,1$, a linear operator
$T$ is a bounded map from the Banach space $X_{j}$ to the Banach
space $Y_{j}$ then, for each $\theta$, $0<\theta<1$, $T$ is a
bounded map between the complex interpolation spaces $\left[X_{0},X_{1}\right]_{\theta}$
and $\left[Y_{0},Y_{1}\right]_{\theta}$. 

Alberto Calder\'on, in his foundational presentation of this material
fifty-one years ago \cite{CalderonA1964}, also proved the following
companion result: 

\textbf{Compactness Theorem/Question:} Furthermore \textbf{\textit{in
some cases}}, if $T$ is also a compact map from $X_{0}$ to $Y_{0}$,
then, for each $\theta$, $T$ is a compact map from $\left[X_{0},X_{1}\right]_{\theta}$
to $\left[Y_{0},Y_{1}\right]_{\theta}$. 

The fundamental question of exactly which cases could be covered by
such a result was not resolved then, and is still not resolved. In
a previous paper \cite{aaCwikelMRochbergR201-A-4arXiv} we surveyed
several of the partial answers which have been obtained to this question,
with particular emphasis on the work of Nigel Kalton in a joint paper
\cite{aaaCwKa} with one of us. This is a preliminary version of a
set of lecture notes which will be a sequel to \cite{aaCwikelMRochbergR201-A-4arXiv}.
In them, for the most part, we will amplify upon various technical
details of the contents of \cite{aaaCwKa}. For example we plan to
give a more explicit explanation of why the positive answer in \cite{aaaCwKa}
to the above question when $\left(X_{0},X_{1}\right)$ is a couple
of lattices holds without any restriction on those lattices, and we
also plan to provide more detailed versions of some of the other proofs
in that paper. The main purpose of this preliminary version is to
present two apparently new small results, pointing out a previously
unnoticed particular case where the answer to the above question is
affirmative. As our title suggests, this and future versions of these
notes are intended to be more accessible to graduate students than
a usual research article. 
\end{abstract}
\maketitle

\section{\label{sec:IntroB}Introduction}

This is a preliminary version of a set of lecture notes which will
be a mostly technical sequel to \cite{aaCwikelMRochbergR201-A-4arXiv},
only intended for readers who are familiar with the contents of \cite{aaCwikelMRochbergR201-A-4arXiv},
or at least with Sections 1 to 3 of that paper. We shall unhesitatingly
use any and all of the notation, terminology and conventions introduced
there, usually without further explanation. This means, among other
things, that all Banach spaces considered here will be over the complex
field. It will also be necessary to consult a number of other references.
In particular, we will assume that the reader is familiar with the
definitions and various properties of the complex interpolation spaces
$\left[A_{0},A_{1}\right]_{\theta}$ and $\left[A_{0},A_{1}\right]^{\theta}$
and the auxiliary spaces $\mathcal{F}\left(A_{0},A_{1}\right)$ and
$\overline{\mathcal{F}}\left(A_{0},A_{1}\right)$ of $A_{0}+A_{1}$
valued functions on the strip $\mathbb{S}=\left\{ z\in\mathbb{C}:0\le\mathrm{Re}\, z\le1\right\} $
which are used in their construction, as can be found in the earlier
sections of \cite{CalderonA1964} or (with slightly different notation)
in Chapter 4 of \cite{BerghLofstrom}. 

Our main goal will be to provide some extra details about some tools
which might ultimately be relevant for helping to answer a 51 year
old question, which, as in \cite{aaCwikelMRochbergR201-A-4arXiv},
we will refer to as \textit{Question CIC.} So far, those who have
considered this question have had to content themselves with finding
various special cases in which its answer is affirmative. I.e., rather
than determining whether or not $\left(\ast,\ast\right)\blacktriangleright(\ast,\ast)$
holds, they have merely found various examples, sometimes quite large
families of examples, of complex Banach couples $\left(X_{0},X_{1}\right)$
and $\left(Y_{0},Y_{1}\right)$ which satisfy $\left(X_{0},X_{1}\right)\blacktriangleright\left(Y_{0},Y_{1}\right)$. 

We have deliberately used rather exotic notation in the preceding
paragraph, in order to further emphasize that we assume familiarity
with \cite{aaCwikelMRochbergR201-A-4arXiv}, where that notation is
explained, and where some history of this question is recalled. 

Given the fact that Question CIC has remained open for five decades,
together with the fact that in the nineteen years since the publication
of \cite{aaaCwKa} there has been rather little further progress towards
answering it, we feel that there is a case for carefully looking again
at the details of \cite{aaaCwKa}. But we have deferred doing this
until a later version of these notes. 

Our own recent rereading of \cite{aaaCwKa} has born some modest fruit,
by prompting us to discover at least one new family of Banach couples
$\left(B_{0},B_{1}\right)$ for which $\left(\ast,\ast\right)\blacktriangleright\left(B_{0},B_{1}\right)$.
Since we are impatient to report some progress, even if rather small,
in the battle with Question CIC, we have prepared this preliminary
version of our notes mainly for this purpose.

Although we have by and large adopted a format which we hope will
more conveniently accessible to and useful for graduate students,
we should hasten to add that we most certainly do \textbf{\textit{not}}
wish to encourage anyone beginning a mathematical career to choose
answering Question CIC as the main topic for her or his dissertation.
That would really be a ``high risk trajectory''.

We originally wrote \cite{aaCwikelMRochbergR201-A-4arXiv} as the
fourth of a series of papers intended to describe some small part
of the most impressive body of mathematical research created by our
brilliant colleague and dear friend Nigel Kalton. The first three
of these papers were coauthored with Mario Milman. All of them have
been posted on the arXiv. We have submitted material taken from those
papers to also possibly appear in a ``Selecta'' volume to be published
in Nigel's memory.

\section{\label{sec:ViaSchauder}A new result via Schauder's theorem}

In this section, in Corollary \ref{cor:NewUMD}, we find, as promised,
one more class of Banach couples $\left(B_{0},B{}_{1}\right)$ for
which $\left(*,*\right)\blacktriangleright\left(B_{0},B_{1}\right)$,
i.e., for which $\left(A_{0},A_{1}\right)\blacktriangleright\left(B_{0},B_{1}\right)$
for all Banach couples $\left(A_{0},A_{1}\right)$. The result is
not entirely surprising, in view of Theorem 9 of \cite[p.~271]{aaaCwKa}
and Schauder's classical theorem about adjoints of compact operators.
Indeed our result emerges as a consequence (see Theorem \ref{thm:DualizeQ1})
of a Schauder-like theorem for complex interpolation of compact operators. 

It will be convenient (as in \cite[p.~22]{CwikelIntoLattices}) to
use the notation $T:\left(A_{0},A_{1}\right)\overset{c,b}{\to}\left(B_{0},B_{1}\right)$,
to mean that $\left(A_{0},A_{1}\right)$ and $\left(B_{0},B_{1}\right)$
are Banach couples and that $T:A_{0}+A_{1}\to B_{0}+B_{1}$ is a linear
operator which maps $A_{0}$ compactly into $B_{0}$ and $A_{1}$
boundedly into $B_{1}$. 

To begin our discussion we take the liberty of recalling an obvious
fact which has been frequently used, in \cite{aaaCwKa} and elsewhere: 
\begin{fact}
\label{fact:Clintersect}Let $\left(A_{0},A_{1}\right)$ and $\left(B_{0},B_{1}\right)$
be Banach couples, and let $\left(X_{0},X_{1}\right)$ and $\left(Y_{0},Y_{1}\right)$
be the regular Banach couples obtained by letting $X_{j}$ be the
closure of $A_{0}\cap A_{1}$ in $A_{j}$ and $Y_{j}$ be the closure
of $B_{j}$ in $B_{0}\cap B_{1}$ for $j=0,1$. Then $\left(X_{0},X_{1}\right)\blacktriangleright\left(Y_{0},Y_{1}\right)$
implies that $\left(A_{0},A_{1}\right)\blacktriangleright\left(B_{0},B_{1}\right)$.
\end{fact}
We also take the liberty of recalling that Fact \ref{fact:Clintersect}
is an immediate consequence of the obvious fact that any linear operator
$T$ satisfying $T:\left(A_{0},A_{1}\right)\overset{c,b}{\to}\left(B_{0},B_{1}\right)$
also satisfies $T:\left(X_{0},X_{1}\right)\overset{c,b}{\to}\left(Y_{0},Y_{1}\right)$
and also the very well known fact (see Section 9.3 of \cite[p.~116]{CalderonA1964})
that $\left[A_{0},A_{1}\right]_{\theta}=\left[X_{0},X_{1}\right]_{\theta}$
isometrically (and analogously $\left[B_{0},B_{1}\right]_{\theta}=\left[Y_{0},Y_{1}\right]_{\theta}$
isometrically). (Incidentally, since any linear operator satisfying
$T:\left(X_{0},X_{1}\right)\overset{c,b}{\to}\left(Y_{0},Y_{1}\right)$
can be uniquely extended to be an operator $\widetilde{T}$ satisfying
$\widetilde{T}:\left(A_{0},A_{1}\right)\overset{c,b}{\to}\left(B_{0},B_{1}\right)$,
we can also observe that the two properties $\left(X_{0},X_{1}\right)\blacktriangleright\left(Y_{0},Y_{1}\right)$
and $\left(A_{0},A_{1}\right)\blacktriangleright\left(B_{0},B_{1}\right)$
are in fact equivalent.)

Our discussion continues with the following auxiliary result:
\begin{thm}
\label{thm:Lower2Upper}Let $\left(A_{0},A_{1}\right)$ and $\left(B_{0},B_{1}\right)$
be arbitrary Banach couples. Let $U$ be a bounded linear operator
from $A_{0}+A_{1}$ into $B_{0}+B_{1}$. Suppose that, for some $\theta\in(0,1)$,
the operator $U$ maps $\left[A_{0},A_{1}\right]_{\theta}$ compactly
into $\left[B_{0},B_{1}\right]_{\theta}$. Then $U$ also maps $\left[A_{0},A_{1}\right]^{\theta}$
compactly into $\left[B_{0},B_{1}\right]_{\theta}$. \end{thm}
\begin{rem}
The boundedness of $U:A_{0}+A_{1}\to B_{0}+B_{1}$ will of course
hold whenever the much stronger property $U:\left(A_{0},A_{1}\right)\overset{c,b}{\to}\left(B_{0},B_{1}\right)$
is assumed to hold. We also remark that the conclusion of the theorem
also holds if $\left[A_{0},A_{1}\right]^{\theta}$ is replaced by
the (possibly larger?) Gagliardo closure%
\footnote{This is the space whose unit ball is the closure of the unit ball
of $\left[A_{0},A_{1}\right]_{\theta}$ with respect to the norm of
$A_{0}+A_{1}$.%
} in $A_{0}+A_{1}$ of $\left[A_{0},A_{1}\right]^{\theta}$. The following
proof has some similarities with the proof of (b) in the second paragraph
of the proof of Corollary 7 of \cite[p.~270]{aaaCwKa}. 
\end{rem}
\noindent \textit{Proof.} Let $a$ be an arbitrary element of $\left[A_{0},A_{1}\right]^{\theta}$.
There exists a function $f\in\overline{\mathcal{F}}\left(A_{0},A_{1}\right)$
for which $a=f^{\prime}(\theta)$. Since $f$ is an analytic $A_{0}+A_{1}$
valued function in the open strip $\mathbb{S}^{\circ}$, this means
that $\lim_{h\to0}\ell\left(\frac{f\left(\theta+h\right)-f\left(\theta\right)}{h}-f^{\prime}(\theta)\right)=0$
for every bounded linear functional $\ell$ on $A_{0}+A_{1}$. In
turn this implies that 
\begin{equation}
\lim_{h\to0}\lambda\left(U\left(\frac{f\left(\theta+h\right)-f\left(\theta\right)}{h}\right)-Uf^{\prime}(\theta)\right)=0\label{eq:WkCgce}
\end{equation}
for every bounded linear functional $\lambda$ on $B_{0}+B_{1}$.
It is obvious from the definition%
\footnote{Note that in \cite{BerghLofstrom} the notation $\mathcal{G}\left(A_{0},A_{1}\right)$
is used to denote the space $\overline{\mathcal{F}}\left(A_{0},A_{1}\right)$
of \cite{CalderonA1964}. The same notation $\mathcal{G}\left(A_{0},A_{1}\right)$
has a quite different meaning in \cite{CalderonA1964}.%
} of $\overline{\mathcal{F}}\left(A_{0},A_{1}\right)$ (as already
observed and used long ago in \cite[p.~136]{CalderonA1964} and also
used elsewhere, e.g., in \cite[p.~1006]{CwikelM1978-complex} and
\cite{CwikelDualityLectures}) that, for each $n\in\mathbb{N}$, the
function $f_{n}:\mathbb{S}\to A_{0}+A_{1}$ defined by $f_{n}(z)=ine^{z^{2}/n}\left(f(z+1/in)-f(z)\right)$
is an element of $\mathcal{F}\left(A_{0},A_{1}\right)$ and satisfies
$\left\Vert f_{n}\right\Vert _{\mathcal{F}\left(A_{0},A_{1}\right)}\le e^{1/n}\left\Vert f\right\Vert _{\overline{\mathcal{F}}\left(A_{0},A_{1}\right)}$.
It follows that $\left\{ f_{n}(\theta)\right\} _{n\in\mathbb{N}}$
is a bounded sequence in $\left[A_{0},A_{1}\right]_{\theta}$ and
therefore $\left\{ Uf_{n}(\theta)\right\} _{n\in\mathbb{N}}$ has
a subsequence which converges in $\left[B_{0},B_{1}\right]_{\theta}$
norm to an element $b\in\left[B_{0},B_{1}\right]_{\theta}$. This
convergence must also occur with respect to the norm of $B_{0}+B_{1}$
(since of course $\left[B_{0},B_{1}\right]_{\theta}$ is continuously
embedded in $B_{0}+B_{1}$), and therefore it also occurs with respect
to the weak topology of $B_{0}+B_{1}$. Since the limit of a weakly
convergent sequence is unique, we deduce, using (\ref{eq:WkCgce}),
that $b$ must equal $Uf'(\theta)=Ua$. 

Now let $\left\{ a_{k}\right\} _{k\in\mathbb{N}}$ be an arbitrary
bounded sequence in $\left[A_{0},A_{1}\right]^{\theta}$. The arguments
of the preceding paragraph imply, for each $k\in\mathbb{N}$, that
$Ua_{k}\in\left[B_{0},B_{1}\right]_{\theta}$ and also there that
exists an element $x_{k}\in\left[A_{0},A_{1}\right]_{\theta}$ such
that $\left\Vert x_{k}\right\Vert _{\left[A_{0},A_{1}\right]_{\theta}}\le e\left\Vert a_{k}\right\Vert _{\left[A_{0},A_{1}\right]^{\theta}}$
and $\left\Vert Ua_{k}-Ux_{k}\right\Vert _{\left[B_{0},B_{1}\right]_{\theta}}\le1/k$.
The fact that $U:\left[A_{0},A_{1}\right]_{\theta}\to\left[B_{0},B_{1}\right]_{\theta}$
is compact ensures that $\left\{ Ux_{k}\right\} _{k\in\mathbb{N}}$
has a convergent subsequence in $\left[B_{0},B_{1}\right]_{\theta}$
and therefore that the same is true of $\left\{ Ua_{k}\right\} _{k\in\mathbb{N}}$.
$\qed$ 

Now we can state the main results of this section.
\begin{thm}
\label{thm:DualizeQ1}Let $\left(A_{0},A_{1}\right)$ and $\left(B_{0},B_{1}\right)$
be Banach couples, and let $\left(X_{0},X_{1}\right)$ and $\left(Y_{0},Y_{1}\right)$
be the regular Banach couples obtained by letting $X_{j}$ be the
closure of $A_{0}\cap A_{1}$ in $A_{j}$ and $Y_{j}$ be the closure
of $B_{j}$ in $B_{0}\cap B_{1}$ for $j=0,1$. Suppose that the dual
couples $\left(X_{0}^{*},X_{1}^{*}\right)$ and $\left(Y_{0}^{*},Y_{1}^{*}\right)$
satisfy $\left(Y_{0}^{*},Y_{1}^{*}\right)\blacktriangleright\left(X_{0}^{*},X_{1}^{*}\right)$.
Then $\left(A_{0},A_{1}\right)\blacktriangleright\left(B_{0},B_{1}\right)$.
\end{thm}

\begin{cor}
\label{cor:NewUMD}Let $\left(A_{0},A_{1}\right)$ and $\left(B_{0},B_{1}\right)$
be Banach couples, such that $B_{0}$ is a UMD-space. Then $\left(A_{0},A_{1}\right)\blacktriangleright\left(B_{0},B_{1}\right)$.\end{cor}
\begin{rem}
\label{rem:PreAdjoint} It seems natural to conjecture that a sort
of converse to Theorem \ref{thm:DualizeQ1} also holds, namely that
$\left(A_{0},A_{1}\right)\blacktriangleright\left(B_{0},B_{1}\right)$
or, equivalently, $\left(X_{0},X_{1}\right)\blacktriangleright\left(Y_{0},Y_{1}\right)$
is sufficient to imply that $\left(Y_{0}^{*},Y_{1}^{*}\right)\blacktriangleright\left(X_{0}^{*},X_{1}^{*}\right)$.
However, as explained in the second version of \cite{aaCwikelMRochbergR201-A-4arXiv},
proving this conjecture, even for just one special particular choice
of $\left(X_{0},X_{1}\right)$ and $\left(Y_{0},Y_{1}\right)$, namely
$\left(X_{0},X_{1}\right)=\left(\ell^{1}\left(FL^{1}\right),\ell^{1}\left(FL_{1}^{1}\right)\right)$
and $\left(Y_{0},Y_{1}\right)=\ell^{\infty}\left(FL^{\infty}\right),\ell^{\infty}\left(FL_{1}^{\infty}\right)$,
would be equivalent to obtaining a positive answer for Question CIC.
\end{rem}

\begin{rem}
\label{rem:DualOfUMD}In the proof of Corollary \ref{cor:NewUMD}
we will need to use two fairly obvious facts, namely that the dual
$Y^{*}$ of a UMD-space $Y$ and also every closed subspace of $Y$
are also both UMD-spaces. These are mentioned without proof in several
papers (including \cite{aaaCwKa}). Formal statements of these facts
can be found as parts (v) and (viii) respectively of Theorem 4.5.2
of \cite[p.~145]{AmannH1995} and the proof of (v) there is provided
as one of the consequences of Theorem 4.3.6 on p.~139 of the same
book. It seems that there are easier settings in which to write such
proofs, for example in the context of Corollary 2.18 of \cite[p.~495]{BerksonEGillespieTMuhlyP1986},
if one uses the characterization of UMD-spaces which appears there.
In the future full version of these notes we may perhaps offer a fairly
``self contained'' treatment of the material that we require about
UMD-spaces, working merely in terms of trigonmetric polynomials and
thus bypassing the need for various technical details.
\end{rem}
\smallskip{}

\noindent \textit{Proof of Corollary \ref{cor:NewUMD}. }Let $\left(X_{0},X_{1}\right)$
and $\left(Y_{0},Y_{1}\right)$ be the regular Banach couples obtained
from $\left(A_{0},A_{1}\right)$ and $\left(B_{0},B_{1}\right)$ as
in the statement of Theorem \ref{thm:DualizeQ1}. Since $B_{0}$ is
a UMD-space, so is its closed subspace $Y_{0}$ and therefore so is
the dual $Y_{0}^{*}$ of $Y_{0}$.\textbf{\textit{ }}Accordingly,
Theorem 9 of \cite[p.~271]{aaaCwKa} implies that $\left(Y_{0}^{*},Y_{1}^{*}\right)\blacktriangleright\left(X_{0}^{*},X_{1}^{*}\right)$.
The desired conclusion now follows from Theorem \ref{thm:DualizeQ1}.
$\qed$ 

\smallskip{}

\noindent \textit{Proof of Theorem \ref{thm:DualizeQ1}.} We will
present this proof using the somewhat pedantic language and notation
of \cite{CwikelDualityLectures}. The main ingredient of the proof
is Schauder's theorem about the compactness of the adjoint of a linear
operator. We will use the variant of that theorem which appears as
Theorem 2.7 of \cite[p.~21]{CwikelIntoLattices}. In preparation for
that, as in Section 1 of \cite{CwikelDualityLectures}, we shall introduce
the bilinear functional $\left\langle \cdot,\cdot\right\rangle $
defined on $(X_{0}\cap X_{1})\times(X_{0}\cap X_{1})^{*}$ which defines
the duality between $X_{0}\cap X_{1}$ and its dual. Then, for each
regular intermediate space $X$ with respect to the couple $\left(X_{0},X_{1}\right)$,
we define the space $X^{\#}$ as in Definition 1.4 of \cite[p.~3]{CwikelDualityLectures}.
The Hahn-Banach Theorem, together with the fact (\cite[Fact 1.6, p.~4]{CwikelDualityLectures})
that $X^{\#}$ identifies isometrically with the dual $X^{*}$ of
$X$ imply that 
\begin{equation}
\left\Vert x\right\Vert _{X}=\sup\left\{ \left|\left\langle x,z\right\rangle \right|:z\in X^{\#},\,\left\Vert z\right\Vert _{X^{\#}}\le1\right\} \mbox{ for all }x\in X_{0}\cap X_{1}.\label{eq:Norming}
\end{equation}
In particular we will be considering the cases where $X$ is $X_{0}$
or $X_{1}$ or $\left[X_{0},X_{1}\right]_{\theta}$. We will also
be using the fact that 
\begin{equation}
\left(\left[X_{0},X_{1}\right]_{\theta}\right)^{\#}=\left[X_{0}^{\#},X_{1}^{\#}\right]^{\theta}\mbox{ with equality of norms}\label{eq:caldduthm}
\end{equation}
which, since $X_{0}\cap X_{1}$ is dense in $\left[X_{0},X_{1}\right]_{\theta}$,
is just a reformulation of Calder\'on's duality theorem (\cite{CalderonA1964}
Section 12.1 p.\ 121 and Section 32.1 pp.\ 148--156, or \cite[pp.~98-101]{BerghLofstrom},
or see an alternative proof in Section 2 of \cite{CwikelDualityLectures}).

We will need the exact analogues of these notions for the couple $\left(Y_{0},Y_{1}\right)$.
To avoid confusion we shall use a different notation, namely $\left\langle \!\left\langle \cdot,\cdot\right\rangle \!\right\rangle $
for the bilinear functional defined on $(Y_{0}\cap Y_{1})\times(Y_{0}\cap Y_{1})^{*}$
which defines the duality between $Y_{0}\cap Y_{1}$ and its dual.
For each regular intermediate space $Y$ with respect to the couple
$\left(Y_{0},Y_{1}\right)$ we can unambigously use the notation $Y^{\#}$
for the space defined, as in Definition 1.4 of \cite[p.~3]{CwikelDualityLectures}
(but of course using $\left\langle \!\left\langle \cdot,\cdot\right\rangle \!\right\rangle $).
Here again we will particularly need to consider the cases where $Y$
is $Y_{0}$, $Y_{1}$ or $\left[Y_{0},Y_{1}\right]_{\theta}$. 

Let $T:\left(X_{0},X_{1}\right)\to\left(Y_{0},Y_{1}\right)$ be a
linear operator satisfying $\left\Vert T\right\Vert _{X_{j}\to Y_{j}}\le1$
for $j=0,1$. We shall define the function $h:\left(X_{0}\cap X_{1}\right)\times\left(Y_{0}\cap Y_{1}\right)^{*}\to\mathbb{C}$
by setting 
\begin{equation}
h(a,b)=\left\langle \!\left\langle Ta,b\right\rangle \!\right\rangle \mbox{ for all \ensuremath{a\in X_{0}\cap X_{1}\mbox{ and all }b\in\left(Y_{0}\cap Y_{1}\right)^{*}}}\label{eq:Define-h}
\end{equation}

For each fixed $b\in\left(Y_{0}\cap Y_{1}\right)^{*}$ let $Sb$ be
the linear functional on $X_{0}\cap X_{1}$ which is defined by $Sb(a)=h(a,b)$
for each $a\in X_{0}\cap X_{1}$. Then 
\[
\left|Sb(a)\right|\le\left\Vert T\right\Vert _{X_{0}\cap X_{1}\to Y_{0}\cap Y_{1}}\left\Vert a\right\Vert _{X_{0}\cap X_{1}}\left\Vert b\right\Vert _{\left(Y_{0}\cap Y_{1}\right)^{*}}\le\left\Vert a\right\Vert _{X_{0}\cap X_{1}}\left\Vert b\right\Vert _{\left(Y_{0}\cap Y_{1}\right)^{*}}
\]
which shows that $S:\left(Y_{0}\cap Y_{1}\right)^{*}\to\left(X_{0}\cap X_{1}\right)^{*}$
with $\left\Vert S\right\Vert _{\left(Y_{0}\cap Y_{1}\right)^{*}\to\left(X_{0}\cap X_{1}\right)^{*}}\le1$.
Thus we can rewrite $Sb(a)$ as $\left\langle a,Sb\right\rangle $. 
\begin{rem}
Of course $S$ is essentially the adjoint of $T$. But we need to
apply the operators $T$ and $S$ to several different spaces, i.e.,
it would seem that the bilinear functionals expressing the relevant
dualities for these different spaces are defined on different spaces.
Hence our preference to proceed cautiously, and to use the pedantic
notation of \cite{CwikelDualityLectures}. We should perhaps explicitly
recall that, as remarked after Definition 1.4 of \cite[p.~3]{CwikelDualityLectures},
$\left(X_{0}\cap X_{1}\right)^{*}=\left(X_{0}\cap X_{1}\right)^{\#}$,
and furthermore $\left(X_{0}\cap X_{1}\right)^{\#}=X_{0}^{\#}+X_{1}^{\#}$
(cf.~\cite[Fact 1.9, p.~5]{CwikelDualityLectures}), so that we have
$S:Y_{0}^{\#}+Y_{1}^{\#}\to X_{0}^{\#}+X_{1}^{\#}$. 
\end{rem}
Now let $X$ and $Y$ be regular intermediate spaces with respect
to $\left(X_{0},X_{1}\right)$ and $\left(Y_{0},Y_{1}\right)$ respectively,
and suppose that they are also relative exact interpolation spaces
with respect to these couples. Then $T:X\to Y$ with $\left\Vert T\right\Vert _{X\to Y}\le1$.
It is an easily checked and essentially standard fact that 
\begin{equation}
S:Y^{\#}\to X^{\#}\mbox{ with }\left\Vert S\right\Vert _{Y^{\#}\to X^{\#}}\le1.\label{eq:siioy}
\end{equation}
Let us nevertheless recall the argument which gives (\ref{eq:siioy}).
Recall that $Y^{\#}$ is the subspace of $\left(Y_{0}\cap Y_{1}\right)^{*}$
consisting of those elements $y^{*}\in\left(Y_{0}\cap Y_{1}\right)^{*}$
which satisfy

\[
\left\Vert y^{*}\right\Vert _{Y^{\#}}:=\sup\left\{ \left|\left\langle \!\left\langle y,y^{*}\right\rangle \!\right\rangle \right|:y\in Y_{0}\cap Y_{1},\,\left\Vert y\right\Vert _{Y}\le1\right\} <\infty.
\]

\begin{rem}
Note that, in particular, $\left(Y_{0}\cap Y_{1}\right)^{\#}\overset{1}{=}\left(Y_{0}\cap Y_{1}\right)^{*}$.
In fact what we are doing here (and in \cite{CwikelDualityLectures})
is consistent with the standard procedure used in many papers (and
thus also implicitly in the statement of Theorem \ref{thm:DualizeQ1})
of defining the dual couple $\left(Y_{0}^{*},Y_{1}^{*}\right)$ of
a regular couple $\left(Y_{0},Y_{1}\right)$ by identifying $Y_{0}^{*}$
and $Y_{1}^{*}$ as spaces which are continuously embedded into $\left(Y_{0}\cap Y_{1}\right)^{*}$.
Therefore, if we use that identification, we have $Y_{0}^{*}=Y_{0}^{\#}$
and $Y_{1}^{*}=Y_{1}^{\#}$ with equality of norms. 
\end{rem}
Continuing our verification of (\ref{eq:siioy}), we observe that,
since $Y_{0}\cap Y_{1}$ is continuously embedded in $Y$, it follows
that $Y^{\#}$ is continuously embedded in $\left(Y_{0}\cap Y_{1}\right)^{*}$.
Thus $Sy^{*}$ is a defined element of $\left(X_{0}\cap X_{1}\right)^{*}$
for each $y^{*}\in Y^{\#}$ and, for each $x\in X_{0}\cap X_{1}$,
we have $Tx\in Y_{0}\cap Y_{1}$ and 
\begin{align*}
\left|Sy^{*}(x)\right| & =\left|\left\langle x,Sy^{*}\right\rangle \right|=\left|\left\langle \!\left\langle Tx,y^{*}\right\rangle \!\right\rangle \right|\le\left\Vert Tx\right\Vert _{Y}\left\Vert y^{*}\right\Vert _{Y^{\#}}\\
 & \le\left\Vert x\right\Vert _{X}\left\Vert y^{*}\right\Vert _{Y^{\#}}.
\end{align*}
This means exactly that $Sy^{*}\in X^{\#}$ and $\left\Vert Sy^{*}\right\Vert _{X^{\#}}\le\left\Vert y^{*}\right\Vert _{Y^{\#}}$.
Thus we have established (\ref{eq:siioy}). 

For $X$ and $Y$ having the properties specified above, we now let
$A=X_{0}\cap X_{1}\cap\mathcal{B}_{X}$ and $B=\left(Y_{0}\cap Y_{1}\right)^{*}\cap\mathcal{B}_{Y^{\#}}=\mathcal{B}_{Y^{\#}}$.
(Here, in order to avoid confusion with the set $B$ we are denoting
the closed unit ball of any Banach space $E$ by $\mathcal{B}_{E}$
rather than by $B_{E}$.) Consider the semimetric spaces $\left(A,d_{A}\right)$
and $\left(B,d_{B}\right)$ defined exactly as in Theorem 2.7 of \cite[p.~21]{CwikelIntoLattices},
here using the function $h$ defined, as above, by (\ref{eq:Define-h}),
now restricted of course to $A\times B$. For this choice of $A$,
$B$ and $h$ it follows that, for each $a_{1}$ and $a_{2}$ in $A$,
using the analogue of (\ref{eq:Norming}) for the couple $\left(Y_{0},Y_{1}\right)$,
we have 
\begin{align*}
d_{A}(a_{1},a_{2}) & =\sup\left\{ \left|\left\langle \!\left\langle Ta_{1}-Ta_{2},b\right\rangle \!\right\rangle \right|:b\in\mathcal{B}_{Y^{\#}}\right\} \\
 & =\left\Vert Ta_{1}-Ta_{2}\right\Vert _{Y}.
\end{align*}
Furthermore, for each $b_{1}$ and $b_{2}$ in $B$, in view of the
definition of $X^{\#}$ and of its norm, 
\begin{align*}
d_{B}(b_{1},b_{2}) & =\sup\left\{ \left|\left\langle a,Sb_{1}-Sb_{2}\right\rangle \right|:a\in\mathcal{B}_{X}\cap X_{0}\cap X_{1}\right\} \\
 & =\left\Vert Sb_{1}-Sb_{2}\right\Vert _{X^{\#}}.
\end{align*}
The formula for $d_{B}$ shows that $(B,d_{B})$ is totally bounded
if and only if $S:Y^{\#}\to X^{\#}$ is compact. Since $A$ is dense
in $\mathcal{B}_{X}$, the formula for $d_{A}$ shows that $\left(A,d_{A}\right)$
is totally bounded if and only if $T:X\to Y$ is compact. Thus, Theorem
2.7 of \cite[p.~21]{CwikelIntoLattices} gives us the Schauder type
result that 

\begin{equation}
S:Y^{\#}\to X^{\#}\mbox{ is compact if and only if }T:X\to Y\mbox{ is compact.}\label{eq:Schauder}
\end{equation}

(Let us quickly mention that various quantitative versions of (\ref{eq:Schauder})
expressed in terms of covering numbers, can be obtained by applying
quantitative versions of \cite[Theorem 2.7]{CwikelIntoLattices} which
appear in \cite{CwikelMLevyE2008} or other results referred to on
the first page of \cite{CwikelMLevyE2008}.) 

Now suppose that $T:X_{0}\to Y_{0}$ is compact, in addition to the
other conditions which $T$ was already assumed to satisfy. Then $S:Y_{0}^{\#}\to X_{0}^{\#}$
is compact. The fact that $\left(Y_{0}^{\#},Y_{1}^{\#}\right)\blacktriangleright\left(X_{0}^{\#},X_{1}^{\#}\right)$
then tells us that $S:\left[Y_{0}^{\#},Y_{1}^{\#}\right]_{\theta}\to\left[X_{0}^{\#},X_{1}^{\#}\right]_{\theta}$
is compact. We can now use Theorem \ref{thm:Lower2Upper} together
with the fact (see Section 8 of \cite[p.~116]{CalderonA1964}) that
$\left[X_{0}^{\#},X_{1}^{\#}\right]_{\theta}$ is continuously embedded
into $\left[X_{0}^{\#},X_{1}^{\#}\right]^{\theta}$ to obtain that
$S:\left[Y_{0}^{\#},Y_{1}^{\#}\right]^{\theta}\to\left[X_{0}^{\#},X_{1}^{\#}\right]^{\theta}$
is compact. This, combined with (\ref{eq:caldduthm}) (and of course
its counterpart for $\left(Y_{0},Y_{1}\right)$) and with (\ref{eq:Schauder})
enables us to obtain that $T:\left[X_{0},X_{1}\right]_{\theta}\to\left[Y_{0},Y_{1}\right]_{\theta}$
is compact. This shows that $\left(X_{0},X_{1}\right)\blacktriangleright\left(Y_{0},Y_{1}\right)$
and therefore (cf.~Fact \ref{fact:Clintersect}) completes the proof
of the theorem. $\qed$

\end{document}